\documentclass[10pt]{amsart}
\usepackage{color}
\usepackage[shortlabels]{enumitem}
\usepackage{float}
\usepackage{graphicx}
\usepackage{colortbl}
\usepackage{textcomp}
\usepackage{amssymb}
\usepackage{amsmath}
\usepackage{amsthm}
\usepackage{mathrsfs}
\usepackage{amssymb,amsmath,mathrsfs,latexsym,color,wasysym,graphicx}
\usepackage{tikz}
\usepackage{comment}

\newtheorem{theorem}{Theorem}
\newtheorem{lemma}[theorem]{Lemma}
\newtheorem{corollary}[theorem]{Corollary}
\newtheorem{proposition}[theorem]{Proposition}

\theoremstyle{definition}
\newtheorem{definition}[theorem]{Definition}
\newtheorem{example}[theorem]{Example}
\newtheorem{notation}[theorem]{Notation}
\newtheorem{remark}[theorem]{Remark}
\title[Adjacency, Talented Monoids and Leavitt Path Algebras]{An Adjacency Matrix Perspective of Talented Monoids and Leavitt Path Algebras}

\author{Wolfgang Bock, Alfilgen N. Sebandal}

\AtEndDocument{\bigskip{\footnotesize%
     \noindent  (A.\,Sebandal) \textsc{Mindanao State University - Iligan Institute of Technology, Iligan City, Philippines and Western Sydney University, Sydney, Australia} \par
   \noindent  \textit{E-mail address}:    \noindent\texttt{alfilgen.sebandal@g.msuiit.edu.ph} \par
  \addvspace{\medskipamount}
   \noindent (W.\,Bock) \textsc{Technische Universität Kaiserslautern, \\
   Gottlieb-Daimler-Strasse 48, \\
   67663 Kaiserslautern, Germany} \par
   \noindent\textit{E-mail address}:   \texttt{bock@mathematik.uni-kl.de} \par
}}

\begin{document}
\maketitle
\begin{abstract}
	In this article we establish relationships between Leavitt path algebras, talented monoids and the adjacency matrices of the underlying graphs. We show that indeed the adjacency matrix generates in some sense the group action on the generators of the talented monoid. With the help of this we deduce a form of the aperiodicity index of a graph via the talented monoid. We classify hereditary and saturated subsets via the adjacency matrix. Moreover we give a formula to compute all paths of a given length in a Leavitt path algebra based on the adjacency matrix. In addition we discuss the number of cycles in a graph. In particular we give an equivalent characterization of acylic graphs via the adjacency matrix, the talented monoid and the Leavitt path algebra.	
\end{abstract}

\section{Introduction}
In the early 1960's, W.G. Leavitt \cite{leavitt} defined a class of algebras which is universal with respect to an isomorphism property between finite rank free modules. Corresponding to a positive integer $n$, nowadays, this algebra is known as the Leavitt algebra $L(1,n)$ and  in \cite{A1_1_7}, Leavitt established the simplicity of such algebras. \\
\indent In the paper of G. Abrams and G.A. Pino \cite{AP05}, they constructed an algebra which deals with a directed graph $E$ and a field $K$ and called it the Leavitt path algebra $L_K(E)$. This algebra turned out to be the generalization of $L(1,n)$. In fact, $L(1,n)$ is isomorphic to $L_K(R_n)$ where $R_n$ is the graph having one vertex and $n$ loops at that vertex. This coincides to the simplicity criterion given in \cite{AP05}.\\
\indent To a directed graph $E$, one can naturally associate an abelian monoid $M_E$, called the graph monoid
of $E$. The graph monoids were introduced by Ara, Moreno and Pardo (\cite{LPAbook},\cite{A1_4}) in relations with the
theory of Leavitt path algebras. It was shown that the group completion of $M_E$ is the abelian group generated by isomorphism classes of finitely generated $L_K(E)$-modules $-$ the Grothendieck
group $K_0(L_K(E))$.\\
\indent Another monoid associated to a graph called the talented monoid $T_E$ of $E$, can be considered as the ``time evolution model" of the
monoid $M_E$. In this form it was introduced in \cite{A1} and further studied in \cite{A2, A3_8}. It was shown that the
talented monoid captures certain geometry of the direct graph, and hence the algebraic properties of
the associated Leavitt path algebra. As the group completion of $T_E$ is the graded Grothendieck group
$K^{gr}_0 (L_K(E))$ (see \cite{Ara,A1_6}), the graded classification conjecture \cite{A1_8} can be restated to say that the talented
monoid $T_E$ is a complete invariant for the class of Leavitt path algebras.
In all of these algebraic objects, cycles play a key role, when studying the growth of these objects \cite{zel} or even Lie bracket algebras on such algebras \cite{BS22,Nam}. The geometric structure of the underlying directed graph plays a role in the study of many ring theoretic properties, see e.g.~\cite{Abramsdecade,LPAbook} for an overview. Especially hereditary and saturated sets give rise to the graded ideal structure of the corresponding algebras and the structure of the monoids \cite{A1}.
In this paper we use the adjacency matrix to depict this geometry of the graph and link properties of it to properties of the Leavitt path algebras and talented monoids. Often it is easier to consider matrix relations than the direct description of the objects in the algebras. Note that the use of the adjacency matrix was already used for Leavitt path algebras, see e.g.~\cite{LPAbook}. However a focused study as in this article is unknown to the authors so far.
The paper is organized as follows: In Section 3, we see how the talented monoid is related to the adjacency matrix which is then used to realize the contribution of each vertex in an aperiodic graph in the representation of each of the generators.  Section 4 provides characterization of the hereditary and saturated sets as submatrices of the graph's adjacency matrix. We also see that the $\mathbb{Z}$-series of the talented monoid is represented as a series of submatrices whose lengths coincide. In Section 5, we  describe how cycles can be expressed via the adjacency matrix and give an algorithm how these can be counted. Section 6 is containing a formula using the adjacency matrix for the number of paths in the Leavitt path algebra of certain length.

\section{Preliminaries}

A \emph{directed graph} $E$ is a tuple $(E^0, E^1, r,s)$ where $E^0$ and $E^1$ are sets and $r,s$ are maps from $E^1$ to $E^0$.  The elements of $E^0$ are called $vertices$ and the
elements of $E^1$ $edges$. We think of each $e\in E^1$ as an edge pointing from $s(e)$ to $r(e)$, that is, $s(e)$ is the $source$ of $e$  and $r(e)$ is the $range$ of $e$.  A graph $E$ is \emph{finite} if $E^0$ and $E^1$ are both finite.
If $s^{-1}(v)$ is a finite set for every $v\in E^0$, then the graph is called  $row$-$finite$. A vertex $v$ for which $s^{-1}(v)=\varnothing$ is called a $sink$ and for which $|s^{-1}(v)|< \infty $ is called a finite-emitter. 

A (finite) \emph{path} $p$ in $E$ is a sequence $p=\alpha_1\alpha_2\cdots \alpha_n$ of edges $\alpha_i$ in $E$ such that $r(\alpha_i)=s(\alpha_{i+1})$ for $1\leq i \leq i-1$. If $s(\alpha_1)=r(\alpha_n)$ and $s(e_i)\neq s(e_j)$ for all $i\neq j$, then $p$ is called a $cycle$. Two cycles are said to be $disjoint$ if they do not have any common vertices. An $exit$ for a cycle $C$ is any edge $e\not \in C^1$ with $s(e)\in C^0$. A graph $E$ is said to be a $comet$ graph if $E$ has a unique cycle which has no exit. 

Given a graph, one can associate an algebraic structure that corresponds to the structure offered by the geometry of the graph. One of these structures is the so-called talented monoid,  first introduced and studied in this form by Hazrat and Li \cite{A1}. For a row-finite graph, the  \emph{talented monoid} of $E$, denoted by $T_E$, is the abelian monoid generated by $\lbrace v(i): v\in E^0, i \in \mathbb{Z}\rbrace$, subject to 
$$v(i)= \sum_{e\in s^{-1}(v) }r(e)(i+1)$$
for every $i \in \mathbb{Z}$ and every $v\in E^0$ that is not a sink.

\indent The additive group $\mathbb{Z}$ of integers acts on $T_E $ via monoid automorphisms by shifting indices: For each $n,i\in \mathbb{Z}$ and $v\in E^0$, define $^nv(i)=v(i+n)$, which extends to an action of $\mathbb{Z}$ on $T_E$. Throughout, we will denote elements $v(0)$ in $T_E$ by $v$.
A $\mathbb{Z}$-$order$ ideal of $T_E$ is a submonoid $I$ of $T_E$ that is closed under the action of $\mathbb{Z}$, that is,  $^\alpha x + {^\beta y}\in I$ if and only if $x,y\in I$ for $\alpha, \beta \in \mathbb{Z}$. $I$ is said to be $cyclic$ if for any $x\in I$, there is an $0\neq n\in \mathbb{Z}
$ such that $^nx=x$.

A sequence of $\mathbb{Z}$-order ideals in $T_E$ $$0\subseteq I_1\subseteq I_2\subseteq\cdots \subseteq I_n =T_E$$ is said to be a $\mathbb{Z}$-composition series if each of the quotient monoid $I_i/I_{i-1}$ is a simple $\mathbb{Z}$-order ideal. In \cite{sebandalvilela}, the Jordan-H\"older Theorem for monoids with group action is presented giving justification to the uniqueness up to isomorphism of a composition series, and hence its length. 

The talented monoid has an interesting relationship to Leavitt path algebras. In \cite{LPAbook}, it is shown that graded ideals and so-called $\mathbb{Z}$-order ideals of the talented monoid are both generated by some set of vertices of a row-finite graph. Hence, it is essential to take a look into the graph structure in order to fully understand the monoidic structure of the talented monoid and the algebraic structure of the Leavitt path algebra. In order to do this, various tools can be used, such as the adjacency matrix of the graph.

The \emph{adjacency matrix}  of a finite graph $E$ with $E^0=\{v_1, v_2. \cdots, v_n  \}$, denoted by $Adj(E) \in \text{Mat}(N \times N,\mathbb{N})$ is defined as $$(Adj(E))_{ij} = \left\{ \begin{array}{cc} n & \text{if there are } \, n \, \text{ edges between } \, v_i  \text{ and } v_j\\ 0 &\text{otherwise.} \end{array} \right.$$

Hence, for any $k\in \mathbb{N}$, $(Adj(E)^k)_{ij}$ is the number of paths of length $k$ with source $v_i$ and range $v_j$.

\section{The Talented Monoid and the Adjacency Matrix}

In this section, we look into the geometry of the graph encoded in the adjacency matrix. We derive a relation of the adjacency matrix to the structure of the talented monoid and provide results that realizes the corresponding structures in both concepts.
\begin{remark}
If for a finite graph $E$, $Adj(E)$ has full rank, the graph has no sources and sinks, since every source leads to a zero column vector and every sink leads to a zero row vector in $Adj(E)$.
\end{remark}

The talented monoid has an interesting relation to the adjacency matrix. Particularly we can ``generate'' the group action on the generators of the talented monoid by orders of the adjacency matrix.

\begin{lemma}
Let $E$ be a row-finite graph with $E^0=\{v_1,v_2,\cdots, v_n  \}$, and $Adj(E)$ its adjacency matrix. Then
$$
\left(\begin{array}{c}
     v_1(0) \\
      \vdots \\
      v_n(0)
\end{array}\right) = Adj(E)^k \left(\begin{array}{c}
     v_1(k) \\
      \vdots \\
      v_n(k)
\end{array}\right),\quad \text{for all } k \in \mathbb{N}.
$$
\end{lemma}

\noindent \textbf{Proof:} Let $Adj(E)$ be the the adjacency matrix of a graph $E$. Then for $m \in \mathbb{N}$,
$$
Adj(E) \left(\begin{array}{c}
     v_1(m+1) \\
      \vdots \\
      v_n(m+1)
\end{array}\right) = \left(\begin{array}{c}
     \sum_{k=1}^n \lambda^{(1)}_k v_k(m+1) \\
      \vdots \\
       \sum_{k=1}^n \lambda^{(n)}_k v_k(m+1)
\end{array}\right),
$$
where $\lambda^{(n)}_k$ denotes the number of edges from $v_n$ to $v_k$. But then 
$$
\sum_{k=1}^n \lambda^{(j)}_k v_k(m+1)= \sum_{e\in s^{-1}(v_j) }r(e)(m+1)=v_j(m),
$$
by definition of the talented monoid.
The assertion then follows iteratively. $\hfill \square$

\begin{definition}
 The $outdegree$ of a vertext $v_i$, denoted by $deg(v_i)$, is defined as 
 the number of edges with source $v_i$. We define the $degree$ $matrix$ $Deg(E)$ of $E$ by
$$
Deg(E)= \left(\begin{array}{c  c c} deg(v_1) &  \cdots & 0\\ \vdots
& \ddots & \vdots\\
0 & \cdots & deg(v_n)\end{array}\right).
$$

\end{definition}

\begin{remark}
We can normalize the adjacency matrix by multiplying it from the left with the inverse of the degree matrix
$$
P=Deg(E)^{-1} Adj(E).
$$
The corresponding matrix $P$ is a stochastic matrix, e.g.~the sum of each rows equals one. Stochastic matrices are used to describe transitions in Markov chains, see e.g.~\cite{Asmussen}. Markov chains themselves have a direct relationship to directed graphs. Indeed if the transition matrix is replaced by the adjacency matrix, we obtain a topological Markov chain. For more about topological Markov chains and $C^*$-algebras, see e.g.~\cite{Wagoner}. 
\end{remark}

\begin{definition}
A strongly connected graph is said to be aperiodic if there is no integer $k > 1$ that divides the length of every cycle of the graph. 
\end{definition}
\begin{proposition}
Let $E$ be row-finite graph and $Adj(E)$ its adjacency matrix. Then $E$ is aperiodic if and only if there exists $k_0 \in \mathbb{N}$ such that $Adj(E)^k$ has no zero entries for all $k\geq k_0$.  
\end{proposition}

\noindent \textbf{Proof:} Note first, that by row-finiteness $Adj(E)^k$ has finite entries for any $k\in \mathbb{N}$. It is shown in \cite{Nelson1977} that for a Markov chain we have that if the graph is aperiodic and the Markov chain is transient, there exists a $k_0$ such that $P^k_{i,j}>0$ for all $k\geq k_0$. 
If $Adj(E)_{i,j}=0$ we have immediately $P_{i,j}=0$, since $Deg(E)P=Adj(E)$ and $Deg(E)$ is diagonal. Moreover $P^k_{i,j} =0$ if and only if  $(Deg(E)P)^k_{i,j} =0$. This completes the proof. 
$\hfill \square$

\begin{definition}
The minimal such $k_0\in \mathbb{N}$ is called the $aperiodic$ $index$ of the graph. 
\end{definition}

Based on the above consideration, there exists an up-to-scaling correspondence between Markov chains and the talented monoid.

\begin{remark}

While in Markov chains, due to conservation of probability, the mass distributed from one node has to be conserved, in the corresponding dynamics based on the talented monoid, mass, i.e.~the number of contributors to a description of vertex $v_k(0)$ are not conserved in quantity but rather grow in every step by the factor of the outbound degrees. Indeed we obtain the following relationship:
\end{remark}

\begin{theorem}
Let $k_0$ be the aperiodic index of a strongly connected graph. Then we have $$v(0)=\left(\begin{array}{c}
     v_1(0) \\
      \vdots \\
      v_n(0)
\end{array}\right) = \left(\begin{array}{c}
     \sum_{i=1}^n \lambda^{(1)}_{i,k_0}  v_i(k_0) \\
      \vdots \\
      \sum_{i=1}^n \lambda^{(n)}_{i,k_0} v_{i}(k_0)
\end{array}\right),$$
where $\lambda^{(l)}_{i,k_0}>0$ for all $l$ and $i$.
\end{theorem}

\noindent \textbf{Proof:} 
Since $Adj(E)^{k_0}$ has no zero entries and 
$$
\left(\begin{array}{c}
     v_1(0) \\
      \vdots \\
      v_n(0)
\end{array}\right) = Adj(G)^{k_0} \left(\begin{array}{c}
     v_1(k_0) \\
      \vdots \\
      v_n(k_0)
\end{array}\right),
$$
the assertion follows directly. $\hfill \square$

\begin{corollary}
For a strongly connected graph, the minimal aperiodic index is the minimal number such that every element of the talented monoid can be represented as a linear combination of every other vertex, each contributing nontrivially to the sum. 
\end{corollary}
In other words, for a strongly connected graph, the minimal aperiodic number in the sense of Markov chains, is the least number where each state is ``reached" given the transition is simultaneous and consistent, and is also shown in the talented monoid as the least index that represents the element $v(0)$, which then could be viewed as the state $v$ at time $0$.

\section{Ideals of the Leavitt path algebra and the Adjacency Matrix}

\begin{definition}
\indent For a graph $E$ and a ring $R$ with identity, we define the \textit{Leavitt path algebra} of $E$, denoted by $L_R(E)$, to be the algebra generated by the sets $\{v:v\in E^0\}$, $\{ \alpha :\alpha \in E^1  \}$ and $\{ \alpha^* : \alpha \in E^1  \}$ with coefficients in $R$, subject to the relations
\begin{enumerate}
\item[\textnormal{(V)}]
$v_iv_j=\delta_{i,j}v_i$ for every $v_i, v_j\in E^0$;
\item[\textnormal{(E)}] $s(\alpha)\alpha=\alpha = \alpha r(\alpha)$ and $r(\alpha)\alpha^*=\alpha^*=\alpha^*s(\alpha)$ for all $\alpha \in E^1$;

\item[\textnormal{(CK1)}] $\alpha^*\alpha'=\delta_{\alpha,\alpha'}r(\alpha)$ for all $\alpha, \alpha'\in E^1$;

\item[\textnormal{(CK2)}]$\sum_{ \{\alpha \in E^1 :s(\alpha)=v    \}  } \alpha \alpha^*=v$ for every non-sink, finite-emitter vertex $v$.
\end{enumerate}
\end{definition}
For a graph $E$ and $e\in E^1$, the elements $e^*$ are called the $ghost$ $edges$ where $s(e^*)=r(e)$ and $r(e^*)=s(e)$. The graph obtained from $E$ by adding the ghost edge $e^*$ for every $e\in E^1$ is called the double graph $\hat{E}$ of $E$. In this paper, we are considering Leavitt path algebras over a field $K$.

\begin{definition}
\cite{A1} Let $E$ be a graph. A subset $H\subseteq E^0$ is said to be \emph{hereditary} if for any $e\in E^1$, we have that $s(e)\in H$ implies $r(e)\in H$. A subset $H\subseteq E^0$ is said to be \emph{saturated} if for a regular vertex $v$, $r(s^{-1}(v))\subseteq H$, then $v\in H$. 
\end{definition}

Let $L_K(E)$ be a Leavitt path algebra with coefficients in the field $K$ associated to the row-finite  graph $E$. Denote by $\mathcal{L}^{gr}(L_K(E))$ the lattice of graded ideals of $L_K(E)$, $\mathcal{L}(E)$ the set of hereditary saturated subsets of $E$, and $\mathcal{L}(T_E)$ the lattice of $\mathbb{Z}$-order-ideals of $T_E$. For $H\subseteq E^0$, let $I(H)$ and $\langle H\rangle $  be the graded ideal and the $\mathbb{Z}$-order ideal generated by $H$, respectively.

\begin{theorem}\cite{LPAbook} \label{LPAbook-main}
 Let $E$ be a row-finite graph. Then there is a lattice isomorphism between $\mathcal{L}(E)$ and $\mathcal{L}^{gr}(L_K(E))$ [resp., between $\mathcal{L}(E)$ and $\mathcal{L}(T_E)$, between $\mathcal{L}(T_E)$ and $\mathcal{L}^{gr}(L_F(E))$] given by 
$H\mapsto I(H)$ [resp., $H\mapsto \langle H\rangle $, $\langle H\rangle \mapsto I(H)$]
where $H$ is a hereditary saturated subset of $E$.
\end{theorem}

Therefore, for a row-finite graph, it is enough study the hereditary saturated sets in order to take a look into the $\mathbb{Z}$-order ideals of the talented monoid and the graded ideals of the Leavitt path algebra.

\begin{definition}
For a graph $E$ and $H\subseteq E^0$, the $induced$ $subgraph$ of $H$ , which we shall denote by $E[H]$ is the subgraph of $E$ with
\begin{center}$E[H]^0=H$ and $E[H]^1=\{e\in E^1: s(e), r(e)\in H\}$. 
\end{center}
For convenience, for $H\subseteq E^0$, we shall denote the adjacency matrix of $E[H]$ by $Adj(H)$, that is, $Adj(H):= Adj(E[H])$. 
\end{definition}

\begin{definition}
Let $M$ be an $n\times n$ matrix. An $m\times k$ $submatrix$ of $N$ of $M$, which we shall denote by $N\subseteq M$, is obtained by selecting rows $i_1, i_2, \cdots , i_m$ and columns $j_1,j_2, \cdots , j_k$ of $M$ and forming a matrix using these entries, in the same relative positions, that is, preserving the order of the rows and columns as in $M$.

When the selected rows and columns are $1,2,\cdots,m$ and  $1,2,\cdots, k$, respectively, $N$ is said to be a $formal$ $submatrix$ of $M$, which we denote by $N\subseteq_f M$. In other words, $N$ is an upper left submatrix of $M$, that is, $N_{i,j}=M_{i,j}$ for all $1\leq i\leq m$ and $1\leq j\leq k$. 

An $m\times m$ submatrix $N$ of $M$ is said to be $principal$, which we shall denote by $N\subseteq_p M$, if $N$ is obtained by selecting rows $i_1, i_2, \cdots, i_m$ and the same columns $i_1, i_2, \cdots, i_m$. 

If $N$ is both a formal principal submatrix of $M$, we shall denote this by $N\subseteq_{fp}M$.

\end{definition}

\begin{example} Consider the matrices $M$, $N_1$ and $N_2$ 
\begin{center}$
M=\left(\begin{array}{c c c} 
0 & 1 & 0  \\
0 & 0 & 1  \\
1& 0 & 0  \\
\end{array}
\right)$, $N_1=\left(\begin{array}{c c} 
0 & 1  \\
0 & 0   \\
\end{array}
\right)$ and $N_2=\left(\begin{array}{c  c} 
1 & 0  \\
0 & 0  \\
\end{array}
\right)$.
\end{center}
Then $N_1\subseteq_{fp} M$ by selecting rows and columns $1$ and $2$, and $N_2\subseteq M$ by selecting rows $1$ and $3$, and columns $2$ and $3$. In Section $5$, we shall denote the selection of rows and columns formally. 
\end{example}

For a graph $E$ and $H\subseteq E^0,$ it directly follows that $Adj(H) \subseteq_p Adj(E)$  by the definition of $E[H]$. We then have the following characterization of a hereditary saturated set in the framework of adjacency matrices.

\begin{theorem} \label{theoremz1}
Let $E$ be a graph and $H\subseteq V$. Then $H$ is a hereditary saturated set if and only if there exists a permutation on $E^0$ such that the adjacency matrix of $E$ could be written of the form 
$$ Adj(E)=\left(\begin{array}{c c} Adj(H) &0 \\
A& B
\end{array}
\right),$$ where for each $i$, $A_{i,s}=0$ for all $s$ if $B_{i,t}=0$ for all $t$.
\end{theorem}

\noindent \textbf{Proof:}
Let $E^0=\{v_1,v_2,\cdots , v_n\}$ where $H=\{v_1,v_2,\cdots, v_m\}$, $m\leq n$.

Suppose $H$ is a hereditary saturated set. By how we set the elements of $E^0$, it follows that the entry $Adj(E)_{i,j}$ where $i,j\leq m$ corresponds to vertices $v_i,v_j\in H$. Then $Adj(H)$ is an $m\times m$ matrix and 
$$ Adj(E)=\left(\begin{array}{c c} Adj(H) &C \\
A& B
\end{array}
\right),$$
where $A$ is $(n-m)\times m$, $B$ is $(n-m)  \times (n-m)$, and $C$ is $m\times (n-m)$.

For convenience, we set
\begin{center}
    $Adj(E)_{i,j}=C_{i,j}$,~ $Adj(E)_{s,t}=A_{s,t}$~ and~ $Adj(E)_{x,y}=B_{x,y}$
\end{center}
for $i,t=1,2,\cdots, m$ and $j,s,x,y=m+1,\cdots , n$.

Suppose $C_{i,j}\neq 0$ for some $i=1,\cdots,m$ and $j=m+1,\cdots , n$. Then there is an edge $e$ with $s(e)=v_i$ and $r(e)=v_j$. However, $v_i\in H$ since $i\leq m$. By hereditariness of $H$, it follows that $v_j\in H$. This is a contradiction since $j>m$ which implies $v_j\not \in H$. Hence, $C_{i,j}=0$ for all such $i,j$. 

Now,let $x=m+1, \cdots , n$ and suppose $B_{x,y}=0$ for all $y=m+1, \cdots ,n$.  Then for all $x, y, =m+1, \cdots , n$, $v_x$ and $v_y$ are not connected by an edge. In contrary, suppose there exists $t_0=1,2, \cdots, m$ such that  $A_{x,t_0}\geq 0$. Then there exists an edge $d$ such that $s(d)=v_x$ and $r(d)=v_{t_0}$. Note that since $t_0\leq m$, $v_{t_0}\in H$.

Let $H'=\{v_t:A_{x,t}\geq 0\}$. Then $v_{t_0}\in H'$ which implies $H'\neq \varnothing$. Also, we have $H'=\{v_t: \textnormal{there~is~} d\in E^0 \textnormal{~with~} s(d)=v_x \textnormal{~and~} r(d)=v_t \}\neq \varnothing$ and $r(s^{-1}(v_x))=H'\subseteq H$, since for all such $t$ with $v_t\in H'$, $t\leq m$ which implies $v_t\in H$. Since $H$ is saturated, it follows that $v_x\in H$, a contradiction. Thus, $A_{x,t}=0$ for all $t=1,\cdots, m$. 

Conversely, suppose there exists a permutation on $E^0$ such that the adjacency matrix of $E$ could be written of the form 
$$ Adj(E)=\left(\begin{array}{c c} Adj(H) &0 \\
A& B
\end{array}
\right),$$ where for each $i$,  $A_{i,s}=0$ for all $s$ if $ B_{i,t}=0$ for all $t$.

Let $v_i\in H$ and $v_j\in E^0$ such that there exists an edge $d$ with $s(d)=v_i$ and $r(v_j)=d$. Then $i\leq m$ and $Adj(E)_{i,j}\geq 1$. Now, since $Adj(E)_{i,s}=0$ for all $s>m$, and $Adj(E)_{i,t}=Adj(H)_{i,t}$ for all $t\leq m$ it follows that $j\leq m$. Hence, $v_j\in H$. Accordingly, $H$ is hereditary. 

Now, let $$~~~~~~~~\{v_{j_1},v_{j_2}, \cdots v_{j_s}\} = r(s^{-1}(v_i))\subseteq H.~~~~~~~~(*)$$
 $j_t\leq m$ for all $t\leq s$. Suppose $v_i\not\in H$, that is $i>m$.  Then $Adj(E)_{i,j_t}=A_{i,j_t}\geq 1$ and $Adj(E)_{i,x}=B_{i,x}$ for all $x>m$ since $i>m$, $j_t\leq m$.  However, $(*)$ also implies that $B_{i,x}=0$ for all $x>m$. This is a contradiction to the description of $A$ and $B$. Thus, $v_i\in H$. Accordingly, $H$ is saturated. $\hfill \square$
 
\begin{definition}
For hereditary saturated subsets $H_1$ and $H_2$ of a graph $E$ with $H_1\subseteq H_2$, define the $quotient$ $graph$ $H_2/H_1$ as a graph such that $(H_2/H_1)^0=H_2\setminus H_1$ and $(H_2/H_1)^1=\{e\in E^1: s(e)\in H_2, r(e)\not \in H_1  \}$. The source and range maps of $H_2/H_1$ are restricted from the graph $E$.
\end{definition}

Notice that for a graph $E$ and  hereditary saturated set $H$, the matrix $B$ in Theorem \ref{theoremz1} shall correspond to the adjacency matrix of the quotient graph $E/H$. 
 
The following corollaries are proved in Theorem \ref{theoremz1} above. 
 
\begin{corollary}
\label{corollaryz2} Let $E$ be graph and $H\subseteq V$. Then $H$ is a hereditary set if and only if there exists a permutation on $E^0$ such that the adjacency matrix of $E$ could be written of the form 
$$ Adj(E)=\left(\begin{array}{c c} Adj(H) &0 \\
A& B
\end{array}
\right).$$ 

\end{corollary}

\begin{corollary}
\label{corollaryz3} Let $E$ be a graph and $H\subseteq V$. Then $H$ is a  saturated set if and only if there exists a permutation on $E^0$ such that the adjacency matrix of $E$ could be written of the form 
$$ Adj(E)=\left(\begin{array}{c c} Adj(H) &C \\
A& B
\end{array}
\right),$$ where for each $i$, $A_{i,s}=0$ for all $s$ if $B_{i,t}=0$ for all $t$.
\end{corollary}

 Let $E$ be graph. Corollaries \ref{corollaryz2} and \ref{corollaryz3} give rise to the following definitions.

\begin{definition}
A formal principal submatrix $N$ of $M$ is said to be $hereditary$ if 
$$ M=\left(\begin{array}{c c} N &0 \\
A& B
\end{array}
\right),$$
and is said to be $saturated$ if
$$ M=\left(\begin{array}{c c} N &C \\
A& B
\end{array}
\right)$$
where for each $i$, $A_{i,s}=0$ for all $s$ if $B_{i,t}=0$ for all $t$,
\end{definition}
 
Based on the definitions above, we then have the following corollary of Theorem \ref{theoremz1}.

\begin{corollary}\label{corollaryx3}
Let $E$ be a graph and $H\subseteq E^0$. Then $H$ is hereditary [respectively, saturated] if and only if for a suitable permutation of elements of $ E^0$, $Adj(H)$ is hereditary [respectively, saturated]. Hence, there is a one-to-one correspondence between the hereditary saturated sets in $E$ and the hereditary saturated  submatrices of $Adj(E)$.
\end{corollary}
\begin{remark} Note that 

$$
Adj(E)^k =\left(\begin{array}{c c} Adj(H)^k &0 \\
K&B^k
\end{array}
\right).
$$
Hence a graph with a hereditary  saturated set $H$ can never be aperiodic. 

Note that especially we have for all $v_i\in H$
$$\left(\begin{array}{c}
     v_1(0) \\
      \vdots \\
      v_m(0)
\end{array}\right) = Adj(H)^k \left(\begin{array}{c}
     v_1(k) \\
      \vdots \\
      v_m(k)
\end{array}\right), $$
for all $k \in \mathbb{N}$.

\end{remark}
  
For a graph $E$, we view vertices $v$ as elements $v(0)$ in its talented monoid $T_E$. Hence, for a $\mathbb{Z}$-order ideal $I$ of $T_E$, by $I\cap E^0$, we shall mean vertices $v$ such that $v(0)\in I$. Now, since every $\mathbb{Z}$-order ideal of the talented monoid is generated by some hereditary saturated sets, the following corollary directly follows.

\begin{corollary}
\label{corollaryy1} Let $E$ be a row-finite graph and $T_E$ its talented monoid. Then  $I\subseteq T_E$ is a $\mathbb{Z}$-order ideal of $T_E$ and only if there exists a permutation on $E^0$ such that $Adj(I\cap E^0)$ is a hereditary saturated submatrix of $Adj(E)$.
\end{corollary}

\begin{definition}
Let $M$ be a matrix. A $matrix$ series for $M$ is a sequence of hereditary and saturated submatices
$$~~~~~0 
\subseteq_{fp} N_1\subseteq_{fp} N_2 \subseteq_{fp} N_3\subseteq_{fp}\cdots  \subseteq_{fp} N_n=M.~~~~~(*)$$
We call a matrix series for $M$ a $matrix$ $composition$ $series$ if for each $i$ with $N_i\subseteq_{fp} N\subseteq_{fp} N_{i+1}$ or for some hereditary and saturated submatix $N$ of $M$, then $N=N_i$ or $N=N_{i+1}$. If $(*)$ is a matrix composition  series, then $n$ is said to be the $length$ of the composition series.
\end{definition}

Let $E$ be a graph and $H_1$ and $H_2$ be two hereditary saturated sets in $E$ with $H_1\subseteq H_2$. Then by applying Theorem \ref{theoremz1} to both $H_1$ and $H_2$, we have
$$Adj(E)=\left(\begin{matrix}\left(\begin{matrix} Adj(H_1) &0 \\
A_2&Adj(H_2/H_1)\end{matrix} \right) &0 \\
A_3&Adj(E/H_2)\end{matrix}\right)$$
with
$$Adj(H_2)= \left(\begin{matrix} Adj(H_1) &0 \\
A_2&Adj(H_2/H_1)\end{matrix} \right).$$

 Let $H_1\subseteq  H_2\subseteq \cdots \subseteq H_n$ be a sequence of hereditary saturated sets of $E$. Then for each $k$, $Adj(H_{k-1})\subseteq_{fp} Adj(H_k)$ and recursively applying Theorem \ref{theoremz1} for each subgraph in the sequence, the following corollary directly follows.

 \begin{corollary} \label{corollaryz4}
 Let $E$ be a graph and $H_1\subseteq H_2\subseteq \cdots \subseteq H_n$ be a sequence of hereditary and saturated sets in $E$. Then for each $k=2,3,\cdots,n$,
 $$ Adj(H_k)=\left(\begin{array}{c c} Adj(H_{k-1}) &0 \\
A_k& Adj(H_k/H_{k-1})
\end{array}
\right),$$ where for each $k$ and $i$, $(A_k)_{i,s}=0$ for all $s$ if $ Adj(H_k/H_{k-1})_{i,t}=0$ for all $t$. More precisely up to suitable permutation, 

 $$ Adj(E)=\left(\begin{matrix} \left(
\begin{matrix}
    \left(\begin{matrix}  \left(\begin{matrix}\left(\begin{matrix} Adj(H_1) &0 \\
A_2&B_2\end{matrix} \right) &0 \\
A_3&B_3\end{matrix}\right) &0 \\
A_4&B_4\end{matrix}\right) & \cdots & 0\\
 \vdots & \ddots \\
 A_n & & B_n
\end{matrix}\right) &0 \\
A_{n+1}&B_{n+1}\end{matrix} \right)$$
where for every $k$, $B_k=Adj(H_k/H_{k-1})$, $H_{n+1}=E$, and  also for every $i$, $(A_k)_{i,s}=0$ for all $s$ if $(B_k)_{i,t}=0$ for all $t$. 
 \end{corollary}

If for each $i$, no hereditary saturated subset $H$ with $H_i\subsetneq H\subsetneq H_{i+1}$, then the matrix series $$0 
\subseteq_{fp} 
Adj(H_1)\subseteq_{fp} Adj(H_2) \subseteq_{fp} Adj(H_3)\subseteq_{fp} \cdots \subseteq_{fp} Adj(H_n)\subseteq_{fp} Adj(E)$$ is a matrix composition series for 
$Adj(E)$ as  presented in Corollary \ref{corollaryz4}. Now, it follows that for each $i$, $\langle H_i\rangle / \langle H_{i-1}\rangle$ is a simple $\mathbb{Z}$-order ideal of $T_E$. Hence,  this corresponds to a composition series in the talented monoid
$$ 0 
\subseteq \langle H_1\rangle \subseteq \langle H_2\rangle \subseteq \langle H_3\rangle \subseteq \cdots \subseteq \langle H_n\rangle \subseteq T_E.$$
By the uniqueness of the composition series for a monoid with a group action is unique by the Jordan-H\"older Theorem \cite{sebandalvilela} and by Corollary \ref{corollaryx3}, we have the following corollary.

\begin{corollary}
Let $E$ be a graph and $T_E$ its talented monoid. Then there is a one-to-one correspondence between a matrix composition series for $Adj(E)$ and the composition series for $T_E$. Thus, the matrix composition series for $Adj(E)$ is unique and has length equal to the length of the composition series of $T_E$. 
\end{corollary}

\begin{example}
 Consider the following graph $E$ and its adjacency matrix $Adj(E)$.

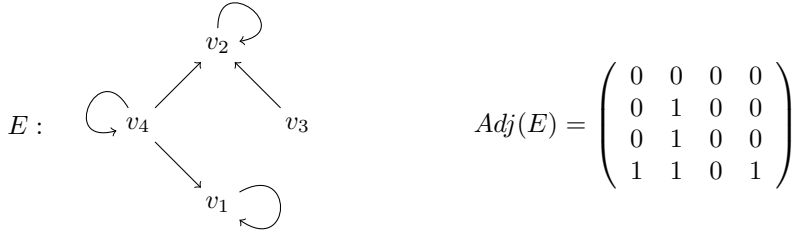
\begin{center}
\begin{tikzpicture}[node distance={15mm}] 
 \node (3) {$v_3$};
\node (2) [above left of=3] {$v_2$};
\node (4) [below left of=2] {$v_4$}; 
\node (1) [below right of=4] {$v_1$}; 
\node (0) [left of =4]{$E:$};
\node (5) [right of =3]{~};
\node (7) [right of =5]{~};
\node (6) [right of =7]{$
Adj(E) =\left(\begin{array}{c c c c } 
0 & 0 & 0  & 0 \\
0 & 1 & 0 & 0 \\
0& 1 & 0  & 0 \\
1 & 1 & 0  & 1 
\end{array}
\right)
$};

\draw [->] (4) to [out=120,in=195,looseness=6] (4);
\draw [->] (1) to [out=30,in=-35,looseness=6] (1);
\draw [->] (2) to [out=90,in=10,looseness=6] (2);
\draw[->] (3) -- (2); 
\draw[->] (4) -- (2); 
\draw[->] (4) -- (1); 
\end{tikzpicture}
\end{center}

Notice that if we have a  sequence formal principal submatrices of $Adj(E)$ given by 
$$   \left ( \begin{matrix} \left ( \begin{matrix}
 \left(\begin{matrix} (0) & 0 \\ 0&1\end{matrix} \right)   & \begin{matrix}
     0\\0
 \end{matrix}\\
 \begin{matrix}
    ~ ~~0&~1
 \end{matrix} & 0~
\end{matrix} \right) & \begin{matrix}
    0\\0\\0
\end{matrix}   \\
\begin{matrix}
    ~~~~1&~1&~~0
\end{matrix} & 1~\end{matrix} \right), $$
this corresponds to the sequence of subsets of $H_1\subseteq H_2 \subseteq  H_3 \subseteq E^0$
where $H_1=\{v_1\}$, $H_2=\{v_1,v_2\}$ and $H_3=\{v_1,v_2,v_3\}$.
However, $H_2$ hereditary but not saturated in $E$. This could be seen in the adjacency matrix of $E$ where the corresponding submatrix of $ \{v_1,v_2\}$ is also hereditary but not a saturated submatrix of $Adj(E)$. In the talented monoid, this corresponds to a sequence of $\mathbb{Z}$-order ideals 
$\langle H_1 \rangle \subseteq \langle H_2 \rangle \subseteq  \langle H_3 \rangle\subseteq T_E$ 
which is also not a composition series for $T_E$ as $ \langle H_2 \rangle=\langle H_3 \rangle $ implying that $ \langle H_2 \rangle/\langle H_3 \rangle =0$. 

It could be seen that the matrix composition series of $Adj(E)$ is given by 
$$\left(\begin{matrix} \left(\begin{matrix} (0) & 0 &0 \\ 0&1&0\\
0&1&0\end{matrix} \right)& \begin{matrix}
    0\\0\\0
\end{matrix}\\
\begin{matrix}
~~~1&~1&0  
\end{matrix} & 1~\end{matrix} \right)  $$
which also corresponds to a composition series in $T_E$
$ 0\subseteq \langle v_1\rangle \subseteq \langle v_1,v_2,v_3 \rangle \subseteq T_E$, both having the same length $3$.  
\end{example}

\section{Matrix representation of cycles}

Cycles are determining the growth of the Leavitt path algebra, which is given in a well-known statement of Alahmedi et.~al.~\cite{zel} on the Gelfand-Kirillov dimension. Particularly it is of interest to determine the number of cycles, a chain of cycles and non-disjoint cycles. In this section, we derive how cycles can be found in the adjacency matrix of a graph. In the same manner, we also provide an algorithm in counting the number of cycles in a graph.

Let $S_n$ denote the permutation group on $n$ letters. A $cyclic$ $permutation$ is any element of $S_n$ of the form $(a_1a_2\cdots a_m)$, $m\leq n$. For $\beta= \{ i_1, i_2, \cdots ,i_m\}$, where $i_j$ are distinct elements in  $I_n=\{1,2,\cdots , n\}$, we denote the set of all cyclic permutations in $S_n$ not fixing $i_j$'s
    by $S_{n,\beta}$. We denote the set of  cyclic permutations of $S_n$ of length $m$ by $S_{n,m}$, that is, $S_{n,m}=\bigcup_{|\beta|=m} S_{n,\beta }  $. We shall relate this to a cycle in a graph.
    
    Let $E$ be a finite graph and  $E^0=\{v_1,v_2,\cdots, v_n\}$. We fix an order of these vertices with $(v_1, v_2, \cdots , v_n)$ and using this order, obtain its adjacency matrix $Adj(E)$ with row and column $i$ corresponding to vertex $v_i$.  We shall call this as an $ordered$ $matrix$ $basis$ for the adjacency matrix $Adj(E)$ described. Denote by  $v_i\rightarrow { v_j}$ if there is an edge $e$ with $s(e)=v_i$ and $r(e)=v_j$. For convenience, a cycle given by vertices $v_i$ with $v_i\rightarrow v_{i+1}$, $v_m\rightarrow v_1$, we shall denote this by 
$v_1\rightarrow \cdots \rightarrow v_m$ instead of $v_1\rightarrow \cdots \rightarrow v_m\rightarrow v_1$.

Consider the graph $E$ below and its adjacency matrix $Adj(E)$.

\begin{center}
\begin{tikzpicture}[node distance={15mm}] 
 \node (3) {$v_3$};
\node (2) [above left of=3] {$v_2$};
\node (1) [below left of=2] {$v_1$}; 
\node (4) [below right of=1] {$v_4$}; 
\node (0) [left of =1]{$E:$};
\node (5) [right of =3]{~};
\node (7) [right of =5]{~};
\node (6) [right of =7]{$
Adj(E) =\left(\begin{array}{c c c  c} 
0 & 1 & 0  & 0 \\
0 & 0 & 1 & 0 \\
0& 0 & 0     & 1\\
1& 0 & 0     & 0
\end{array}
\right)
$};

\draw[->] (1) -- (2);
\draw[->] (2) -- (3); 
\draw[->] (3) -- (4); 
\draw[->] (4) -- (1); 
\end{tikzpicture}
\end{center}

Notice that edge for $ v_i\rightarrow v_j$, corresponds an entry $Adj(E)_{ij}\neq 0$ and we have a cycle $C$ given by $v_1\rightarrow v_2\rightarrow v_3 \rightarrow v_4$. We can view this as $v_1$ ``being sent" to $v_2$, $v_2$ ``being sent" to $v_3$, $v_3$ ``being sent" to $v_4$, and $v_4$ ``being sent" to $v_1$. As a permutation, this shall be $\sigma= (1234)$. For this permutation, we associate the entries $Adj(E)_{12}$, $Adj(E)_{23}$, $Adj(E)_{34}$, and $Adj(E)_{41}$, hence are the entries $Adj(E)_{i, \sigma(i)}$. Multiplying the entries, we get  
$Adj(E)_{12}\cdot Adj(E)_{23}\cdot Adj(E)_{34}\cdot Adj(E)_{41}=1$, which is the number of cycles given by $v_1\rightarrow v_2\rightarrow v_3\rightarrow $. Now, consider the graph $E'$ below where we add edges to $E$ and its adjacency matrix.

\begin{center}
\begin{tikzpicture}[node distance={15mm}] 
 \node (3) {$v_3$};
\node (2) [above left of=3] {$v_2$};
\node (1) [below left of=2] {$v_1$}; 
\node (4) [below right of=1] {$v_4$}; 
\node (0) [left of =1]{$E':$};
\node (5) [right of =3]{~};
\node (7) [right of =5]{~};
\node (6) [right of =7]{$
Adj(E') =\left(\begin{array}{c c c  c} 
0 & 2 & 0  & 0 \\
0 & 0 & 1 & 0 \\
0& 0 & 0     & 1\\
2& 0 & 1     & 0
\end{array}
\right)
$};

\draw[->] (1) -- (2); 
\draw[->] (2) -- (3); 
\draw[->] (3) -- (4); 
\draw[->] (4) -- (1); 

\draw [->] (1) to [out=90,in=180] (2);
\draw[->] (1) -- (2);

\draw [->] (4) to [out=180,in=-90] (1);

\draw [->] (4) to [out=0,in=-90] (3);

\end{tikzpicture}

\end{center}

A cycle given by $v_1\rightarrow v_2\rightarrow v_3 \rightarrow v_4$ again corresponds to permutation $(1234)$ and entries  $Adj(E')_{12}$, $Adj(E')_{23}$, and $Adj(E')_{34}$, and $Adj(E')_{41}$. Multiplying again these entries gives $2\cdot  1 \cdot1\cdot 2\cdot= 4$, the number of cycles given by $v_1\rightarrow v_2\rightarrow v_3 \rightarrow v_4$. 

Consider the permutation $(4321)$ which has the corresponding entries 
$Adj(E')_{43}$, $Adj(E')_{32}$, $Adj(E')_{21}$, and $Adj(E')_{14}$. Multiplying these entries, we get $0$ and we see that indeed, there are no cycles in $E'$ given by $v_4\rightarrow v_3\rightarrow v_2\rightarrow v_1$. We then have the following result. Now, consider the cycle given by $v_3\rightarrow v_4$. Similarly, we can relate this to the permutation $(34)$ of $S_4$ and see that we have $Adj(E')_{34}\cdot Adj(E')_{43}=1$ such cycles. Hence we may correspond all of the cycles given a subset of vertices to a cyclic permutation in $S_4$. 

By the previous visualizations and easy combinatorial arguments on the number of paths given the number of edges, the following result directly follows.

\begin{proposition}
Let $E$ be a finite graph and $(v_1, v_2, \cdots, v_n)$ be its ordered matrix basis for $Adj(E)$. Then the number of cycles given by  $v_{i_1}\rightarrow  v_{i_2}\rightarrow  \cdots \rightarrow  v_{i_m}\rightarrow v_{i_1}$, $v_{i_j}\in E^0$, is  
$$\prod_{i\in  \beta } Adj(E)_{i, \sigma(i)} $$ where $\beta =\{i_1 \cdots i_m\}$ and  $\sigma=( i_1 \cdots i_m)  $. 
Thus, the number of cycles with vertices $v_{i_1}, v_{i_2}, \cdots , v_{i_m}$ is given by 
$$\sum_{\sigma \in S_{n,\beta }}\prod_{i\in  \beta } Adj(E)_{i, \sigma(i)}.$$ 
 Hence, the number of cycles in a graph is given by 
 $$\sum_{\substack{\sigma \in S_{n,\beta } \\ \beta\in P(I_n)}}\prod_{i\in \beta } Adj(E)_{i, \sigma(i)}. $$

\end{proposition}

\begin{example}
    Consider the graph $E$ below and its adjacency matrix $Adj(E)$. 

    \begin{center}
\begin{tikzpicture}[node distance={15mm}] 
 \node (3) {$v_3$};
\node (2) [left of=3] {$v_2$};
\node (1) [left of=2] {$v_1$}; 
\node (0) [left of =1]{$E:$};
\node (5) [right of =3]{~};
\node (7) [right of =5]{~};
\node (6) [right of =7]{$
Adj(E) =\left(\begin{array}{c c c  } 
1 & 1 & 0   \\
1 & 0 & 2  \\
1& 1 & 0  
\end{array}
\right)
$};

\draw[->] (2) -- (1); 
\draw[->] (2) to [out=30,in=150] (3); 

\draw [->] (1) to [out=70,in=110] (2);
\draw [->] (2) to [out=70,in=110] (3);
\draw [->] (3) to [out=230,in=-50] (1);
\draw [->] (3) to [out=200, in= -20](2);
\draw [->] (1) to [out=210,in=-215,looseness=6] (1);

\end{tikzpicture}
\end{center}
Here are the computations of the cycles in $E$ given appropriate $\alpha$ and $\beta$. 
\begin{center}
\begin{tabular}{ c |c c }
 cycles given by & $\displaystyle \prod_{i\in  \beta } Adj(E)_{i, \sigma(i)}$  & \\  \hline
$v_1\rightarrow v_1$  & $Adj(E)_{11}=1$ & \\ 
$v_2\rightarrow v_2$  & $Adj(E)_{22}=0$ & \\
$v_3\rightarrow v_3$  & $Adj(E)_{33}=0$ & \\ 
$v_1\rightarrow v_2$ & $Adj(E)_{12}\cdot Adj(E)_{21}=1$  &  \\  
$v_1\rightarrow v_3$  & $Adj(E)_{13}\cdot Adj(E)_{31}=0$ & \\  
$v_2\rightarrow v_3$  & $Adj(E)_{23}\cdot Adj(E)_{32}=2$ &\\  
$v_1\rightarrow v_2\rightarrow v_3$  & $Adj(E)_{12}\cdot Adj(E)_{23}\cdot Adj(E)_{31}=2$ &\\  
$v_3\rightarrow v_2\rightarrow v_1$  &  $Adj(E)_{32}\cdot Adj(E)_{21}\cdot Adj(E)_{13}=0$& 
\end{tabular}
\end{center}
Hence, there are $6$ cycles in $E$.

\end{example}

Now, the following directly follows from the fact that dimensions of path algebras are completely determined by presence of a cycle in the graph 
{{\cite[Theorem 2.6.17]{LPAbook}}} and that cyclic ideals in the talented monoid are also associated to cycles in the graph 
{\cite[Theorem 3.10]{hazratGK}}.

\begin{proposition}
    Let $E$ be a finite graph and $L_K(E)$, $T_E$, $Adj(E)$ the Leavitt path algebra, talented monoid and adjacency matrix  of $E$, respectively. Then the following are equivalent.
    \begin{enumerate}
        \item $E$ is acyclic.
        \item $L_K(E)$ is a finite-dimensional $K$-algebra.
        \item $E$ has disjoint cycles and $T_E/I$ is not cyclic for all $\mathbb{Z}$-order ideals $I$ . 
        \item If $|E^0|=n$, then for all $\beta \subseteq I_n$ and $\sigma \in S_{n,\beta}$, there exists $i\in \beta$ such that $Adj(E)_{i, \sigma(i)}=0$. 
    \end{enumerate}
 \noindent 
\end{proposition}
\noindent \textbf{Proof:}
Proofs of $(1)\Leftrightarrow (2) \Leftrightarrow (4) $ are straightforward.

\noindent $(1\Leftrightarrow 3)$ Let $E$ be a graph with disjoint cycles and such that $T_E/I$ is not cyclic for any $\mathbb{Z}$-order ideal $I$ of $T_E$. Suppose $E$ has a cycle. Since $E$ is a graph with disjoint cycles, we have a maximal cycle $C$, that is, no other cycles in $E$ is connected to $C$. 

Let $K= \bigcup_{v\in C^0}T(v)$, where $T(v)=\{w\in E^0:w=s(p), v=r(p) \textnormal{~for~some~} p\in Path(E)\}$, that is, all of the vertices in $E$ connected to $C$. Let $H=E^0\setminus K$. We show that $H$ is both hereditary and saturated. 

Suppose $H$ is not hereditary, that is, there exists $x\in H$ and $y\not \in H$ such that $x\rightarrow y$. Then $x\in T(y)$. Since $y\not \in H$, $y\in K$. Hence, $x\in K$, a contradiction. Accordingly, $H$ is hereditary. 

Let $v\in E^0$ such that $s^{-1}(v)\subseteq H$. Suppose $v\not \in H$, that is, $v\in K$. Then $v\in T(w)$ for some $w\in C^0 $. This implies that $s^{-1}(v)\not \subseteq H$, a contradiction. Hence, we must have $v\in H$ and $H$ is saturated. 

Now, we take a look at the quotient graph $E/H$. Clearly, $E/H$ has only one cycle which is the maximal cycle $C$. Now, suppose $C$ has an exit in $E/H$. Then there exists  $e\in E^1$ such that $r(e)\in E/H$. Thus, $r(e)\not \in H=E^0\setminus K$. Hence, $r(e)\in K$, that is, $r(e)\in T(v)$ for some $v\in C^0$. We then a cycle which is not disjoint with $C$, a contradiction. Accordingly, $E/H$ has a unique cycle without an exit, hence a comet. Therefore, $T_{E/H}=T_E/\langle H\rangle $ is a cyclic, a contradiction. Hence, $E$ must not have a cycle.

Conversely, suppose $E$ is acyclic. We are left to show $T_E/I$ is not cyclic for any $\mathbb{Z}$-order ideal $I$. We assume the contrary for some ideal $I_0$. Then $T_E/I=T_{E/H}$ is cyclic where $\langle H\rangle =I$. That is, $E/H$ is a comet graph, which has a cycle. This is a contradiction and the proof is done. \qed

Suppose $E$ is acyclic and $T_E$ has a cyclic  $\mathbb{Z}$-order ideal $I$. By Theorem \ref{LPAbook-main}, $I$ is generated by a hereditary saturated set $H$. Hence, by  Corollary 3.12 of {\cite{hazratGK}},  the full subgraph obtained from $H$ is a comet graph. This implies that $H$, hence also $E$, has a cycle, a contradiction. The converse is proved similarly in the other direction. \qed

We now present adjacency matrix representation of cycles in viewpoint of nilpotency of matrices. 

\begin{theorem}\label{thmcyclicblock}
Let $E$ be a graph with a cycle of length $n$ without an exit. Then there exists a permutation of vertices such that the adjacency matrix of $E$ has the form 
$$
Adj(E) =\left(\begin{array}{c c} N &0 \\
A&B
\end{array}
\right).
$$
where $N$ is a circulant permutation matrix of rank $n$ and for every $m\times m$ submatrix $D$ of $N$ obtained by deleting the $(m+1)^{th},\cdots, n^{th}$ rows and columns of $N$, $D$ is nilpotent of index $m$. 
\end{theorem}

\noindent \textbf{Proof:}
Let $E$ be a graph and $c=e_1e_2\cdots e_n$ a cycle without an exit, of length $n$. Then
$$v_1=s(e_1)=r(e_n),v_2=s(e_2)=r(e_1), \cdots , v_n=s(e_n)=r(e_1).$$
Take the permutation of elements $v_1,v_2,\cdots,v_n,w_{n+1},\cdots, w_k$ of $E^0$ where $w_i\in E^0\setminus\{v_1,v_2,\cdots,v_n\}$. Then the adjacency matrix of $E$ is 
$$
Adj(E) =\left(\begin{array}{c c} N &C \\
A&B
\end{array}
\right)
$$
where the entries of $N$ is a circulant permutation matrix that corresponds to the vertices $v_i$, w.l.o.g.
$$
N =\left(\begin{array}{c c c c c c c} 
0 & 1 & 0 & 0 &\cdots & 0 & 0 \\
0 & 0 & 1 & 0 & \cdots & 0 & 0\\
0 & 0 & 0 & 1& \cdots & 0 & 0\\
0 & 0 & 0 & 0 &\cdots & 0 & 0\\
\vdots & \vdots & \vdots & \vdots & \ddots & \vdots  & \vdots\\
0 & 0 & 0 & 0 & 0 &  0& 1\\
1 & 0 & 0 & 0 & 0 &  0 & 0
\end{array}
\right).
$$
 Then the rank $N$ is $n$. Let $m<n$ and $D=N(1,2,\cdots, m)$.
 
 If $m=1$, then $D=0$ and $D$ is nilpotent. Suppose $m>1$. Then $D$ is of the form 
 $$
D =\left(\begin{array}{c c c c c c c} 
0 & 1 & 0 & 0 &\cdots & 0  & 0\\
0 & 0 & 1 & 0 & \cdots & 0 & 0\\
0 & 0 & 0 & 1& \cdots & 0  & 0\\
0 & 0 & 0 & 0 &\cdots & 0 & 0 \\
\vdots & \vdots & \vdots & \vdots & \ddots & \vdots  & \vdots\\
0 & 0 & 0 & 0 & 0 &  0 & 1\\
0 & 0 & 0 & 0 & 0 &  0 & 0
\end{array}
\right)
$$
and $D^m=0$, that is, $D$ is nilpotent of index $m$.

Suppose there exist $i=1,2,\cdots,n$ and $j=n+1,n+2,\cdots,k$ such that $C_{i,j}\neq 0$, that is, $C_{i,j}=t>0$. Then there exists $e\in E^1$ with $s(e)=v_i$ and $r(e)=w_j$. Thus, $e$ is an exit for the cycle $c$, a contradiction. Hence, $C=0$. $\hfill \square$

\begin{remark}\label{remcyclic}
Notice that in Theorem \ref{thmcyclicblock}, if the cycle has an exit, then we write $$
Adj(E) =\left(\begin{array}{c c} N &C\\
A&B
\end{array}
\right).
$$
where where $C\neq 0$ and $N$ is a circulant permutation matrix of rank $n$ and for every $m\times m$ submatrix $D$ of $N$ obtained by deleting the $(m+1)^{th},\cdots, n^{th}$ rows and columns of $N$, $D$ is nilpotent of index $m$. 
\end{remark}

\begin{remark}
Let $E$ be a graph having a single cycle with an exit. Then there exists a permutation of vertices such that the adjacency matrix of $C$ has the form:
$$
M =\left(\begin{array}{c c} N &0 \\
0&B
\end{array}
\right) +\left(\begin{array}{c c}  0& Out\\
In&0
\end{array}
\right),
$$
where $In$ denote the edges leading to the cycle and $Out$ denotes the exits to vertices outside the cycle. And $N$ is a cirulant permutation matrix.
\end{remark}

\begin{theorem}
\cite{hazratGK} \label{hazratGK-3.10} 
Let $E $ be a finite graph. There is a one-to-one correspondence between
cycles with no exits in $E$ and the cyclic minimal ideals of $T_E$.
\end{theorem}

The following theorem follows from Theorem \ref{hazratGK-3.10} and Corollary \ref{corollaryy1}.

\begin{corollary}
Let $E$ be a graph and $T_E$ its talented monoid. If $T_E$ has a cyclic minimal ideal, then there exists a permutation such that the adjacency matrix of $E$ reads
$$
Adj(E) =\left(\begin{array}{c c} I &0 \\
A&B
\end{array}\right),
$$

with 

$$
I=\left(\begin{array}{c c} N &0 \\
H_1&H_2
\end{array}\right), 
$$
where for each $i$, $A_{i,s}=0$ for all $s$ if $B_{i,t}=0$ for all $t$, $N$ is a circulant permutation matrix of rank $n$ and for every $m\times m$ submatrix $D$ of $N$ obtained by deleting the $(m+1)^{th},\cdots, n^{th}$ rows and columns of $N$, $D$ is nilpotent of index $m$.
\end{corollary}

\noindent \textbf{Proof:} Follows from Theorems \ref{thmcyclicblock}
and  \ref{hazratGK-3.10}, and Corollary \ref{corollaryy1}. $\square$

\section{Number of Paths and Adjacency Matrix}

In this section, we count the number of paths in the double graph $\hat{E}$ of a given length, subject to the Cuntz-Krieger relations (CK1) and (CK2) of the Leavitt Path Algebra. For a graph $E$, the  $opposite$ $graph$ $E^{\text{op}}$ of $E$ is the graph having vertices ${(E^{\text{op}})}^0:=E^0$ and edges ${(E^{\text{op}})}^1:= \{e^*:e\in E^1, s(e)=r(e^*), r(e)=s(e^*)\}$ (also called the $transpose$ $graph$, see \cite{LPAbook}).

\begin{notation}
Let $M\in M_{n\times m}(\mathbb{R})$. For each $s\leq n$ and $t\leq m$, denote 
\begin{center}
    $\displaystyle \| M\| ^r_{s}= \sum_{j=1}^mM_{sj}~~$, $~~\displaystyle\| M\| ^c_{t}= \sum_{i=1}^nM_{it}~~$, ~and $~~\displaystyle\| M\| = \sum_{i,j}M_{ij}$.
 \end{center}   
For a finite graph $E$ and $v_i.v_j\in E^0$, $Adj(E)^k_{ij}$ is the number of paths of length $k$ with source $v_i$ and range $v_j$, for a suitable arrangement of rows and columns. Hence, $\| Adj(E)^k\| ^r_i$ and $\| Adj(E)^k\| ^c_j $ is the number of edges with source $v_i$ and the number of edges with range $v_j$, respectively. Accordingly, for $k=0$, we shall have a path of length $0$ and source and range itself, hence a vertex. Thus, it is imperical to setting  $\| Adj(E)^0\| ^r_i, \| Adj(E)^0\| ^c_j =1$ and $\| Adj(E)^0\| = |E^0| $
for any $i,j$.
\end{notation}

\begin{theorem}\label{pathcomputation}
Let $E$ be a finite graph $A$ its adjacency matrix and $\hat{E}$ its double graph.  Then for each $k>1$, the number of paths  of length $k$  in $L_K(E)$ of the form $\alpha \beta^*$ in $\hat{E}$ is
\begin{eqnarray*} p_k(E) &=& 2\|  A^k\| +\sum_{\substack{s+t=k \\ s,t>0}}\left ( \sum_{j=1}^n\left( \| A^s\| ^c_j \right)  \left( \| A^t\| ^c_j  \right)    \right)\\
&~~~~~~~~~~~~~~~~&~~~~~~~~~~ -\sum_{\substack{s+t=k \\ s,t>0}}\left (\sum_{\| A\| ^r_j=1} \left( \| A^{s-1}\| ^c_j \right)  \left( \| A^{t-1}\| ^c_j  \right)\right )  .\end{eqnarray*} 
\end{theorem}

\noindent \textbf{Proof:}
Let $E$ be a finite graph with $|E^0|=n$. Then then the number of paths of length $0$ and $1$ are $|E^0|=n$ and $2|E^1|$, respectively. We now compute the number of paths of length $k>1$. 

Fix $s,t\in \mathbb{N}$ with $s+t=k$. For each $v_j\in E^0$, consider $\alpha v_j\beta^*$ where $\alpha, \beta\in \text{Path}(E)$ with $r(\alpha)=v_j=s(\beta^*)=r(\beta)$ and $l(\alpha)=s$ and $l(\beta)=t$. Suppose $s=0$ or $t=0$. Then there are exactly $\| A^k\| $ paths $\alpha$ and $\| A^k\| $ paths $\beta^*$. We then obtain $2\| A^k\| $ for such paths. 

Suppose $s\neq 0$ and $t\neq 0$. Note that there are $\| A^s\| _j^c$  paths of length $s$ in $E$ with range $v_j$ and $\| A^t\| _j^c$  paths in $E^{op}$ with source $v_j$. Hence, there are $\| Adj(E)^s\| ^c_j\| Adj(E)^t\| ^c_j$ paths $\alpha \beta^*$ with $r(\alpha)=v_j=s(\beta^*)$ and $l(\alpha)=s$ and $l(\beta)=t$. Taking the sum over all vertices and over all such $s,t$ we obtain 

$$\sum_{\substack{s+t=k \\ s,t>0}}\left ( \sum_{j=1}^n\left( \| A^s\| ^c_j \right)  \left( \| A^t\| ^c_j  \right)    \right).$$

Now, in the reduction of path length by the CK-relations, we only need to consider (CK2) since (CK1) is not utilized in such paths. Notice that (CK2) only contributes to the reduction whenever a vertex emits only $1$ edge. In this case, $v=ee^*$ for some $e\in E^1$. \\
Hence, we look at $\alpha \beta^*$ with $\alpha = d_1d_2\cdots d_{s-1}e$ and $\beta^*=e^*f_1^*f_2^*\cdots f_{t-1}^*$ since
$\alpha \beta^*= d_1d_2\cdots d_{s-1}e e^*f_1^*f_2^*\cdots f_{t-1}^*= d_1d_2\cdots d_{s-1}f_1^*f_2^*\cdots f_{t-1}^*$, a path of length $<k$. We deduct the paths $\alpha v_j\beta^*$ for each $v_j\in E^0$ with $|s^{-1}(v_j)|=1$. This corresponds to paths  $d_1d_2\cdots d_{s-1}e e^*f_1^*f_2^*\cdots f_{t-1}^*= d_1d_2\cdots d_{s-1}f_1^*f_2^*\cdots f_{t-1}^*$ with $r(d_{s-1})=s(f_1^*)=v_j$ having $\| A\| ^r_j=1$. Similar arguments imply that there are

$$\sum_{\substack{s+t=k \\ s,t>0}}\left (\sum_{\| A\| ^r_j=1} \left( \| A^{s-1}\| ^c_j \right)  \left( \| A^{t-1}\| ^c_j  \right)\right )   $$ such paths. 

$\hfill \square$

The elements $pq^*$, with $r(p)=r(q)$ generates the Leavitt path algebra where grading is also determined by such paths. Hence, formulas such as in Theorem \ref{pathcomputation} would be useful in finding the dimension of the subspaces that makes up the entire algebra.

\begin{example}
    We shall calculate the number of path $\alpha\beta^*$ of the following graph using its adjacency matrix shown below. 

    \begin{center}
\begin{tikzpicture}[node distance={15mm}] 
 \node (3) {$v_2$};
\node (2) [left of=3] {$v_1$}; 
\node (0) [left of =2]{$E:$};
\node (7) [right of =3]{~};
\node (8) [right of =7]{~};
\node (6) [right of =8]{$
Adj(E) =\left(\begin{array}{c c  } 
0 & 2   \\
0 & 1     
\end{array}
\right)
$};

\draw [->] (2) to [out=60,in=120] (3);
\draw [->] (2) to [out=-60, in= -120](3);
\draw [->] (3) to [out=35,in=-35,looseness=6] (3);

\end{tikzpicture}
\end{center}
Now, for any $k\geq 1$, we have 
$
Adj(E)^k =\left(\begin{array}{c c  } 
0 & 2   \\
0 & 1     
\end{array}
\right).
$
Hence 
$$\|Adj(E)^k\|=3, \quad \|Adj(E)^k\|^c_1, = 0 \quad \|Adj(E)^k\|^c_2 = 3 $$
for all $k \in \mathbb{N}$.
The number of paths in the $L_K(E)$ of length $k=3$ is hence 
\begin{eqnarray*}
p_3(E)&=&2\|  A^3\| +\sum_{\substack{s+t=3 \\ s,t>0}}\left ( \sum_{j=1}^n\left( \| A^s\| ^c_j \right)  \left( \| A^t\| ^c_j  \right)    \right)\\&&-  \sum_{\substack{s+t=3 \\ s,t>0}}\left (\sum_{\| A\| ^r_j=1} \left( \| A^{s-1}\| ^c_j \right)  \left( \| A^{t-1}\| ^c_j  \right)\right )  \\
&=& 6 +\left (0+0+9+9 \right)-  \left (3+3\right ) =18,
\end{eqnarray*}
since we have the tuples $(s,t)=(1,2)$ and $(s,t)=(2,1)$.
Hence there are 18 paths of lengths $3$, which are of the form $\alpha \beta^*$. 

\end{example}

\section*{Acknowledgements} A.S.~gratefully acknowledges the Ph.D. sandwich scholarship within the DOST-ASTHRDP program by Department of Science and Technology, The Philippines. In addition she thanks Western Sydney University for the hospitality within the sandwich program.


\begin{thebibliography}{20}
\baselineskip=1\baselineskip


\bibitem{Abramsdecade} G. Abrams,\emph{Leavitt path algebras: the first decade.} Bulletin of Mathematical Sciences, 5(1), 59-120, (2015).


\bibitem{AP05} G. Abrams, G., and G. A. Pino,  \emph{The Leavitt path algebra of a graph.} Journal of Algebra, 293(2), 319-334, (2005).

\bibitem{AMP07} G. Abrams,  G. A. Pino, M. S. Molina, \emph{Finite-dimensional Leavitt path algebras.} Journal of Pure and Applied Algebra, 209(3), 753-762. (2007).

\bibitem{abramsmesyan} G. Abrams, Z. Mesyan, \emph{Simple Lie algebras arising from Leavitt path algebras}. Journal of Pure and Applied Algebra, 219, 2302-2313, (2012).

\bibitem{LPAbook} G. Abrams, P. Ara, M. Siles Molina, \emph{Leavitt Path Algebras}, Lecture Notes in Mathematics, vol. 2191, Springer Verlag, 2017.

\bibitem{Alahmedi} A. Alahmedi, H. Alsulami, \emph{On the simplicity of the Lie Algebra of a Leavitt Path Algebra}, Communications in Algebra, vol. 44. (2016).

\bibitem{zel} A. Alahmedi, H. Alsulami, S.K. Jain, E. Zelmanov, \emph{Leavitt path algebras of finite Gelfand-Kirillov dimension},  Journal of Algebra and Its Applications, {\bf 11}, No. 06, 1250225 (2012). 


\bibitem{Ara} P. Ara, and E. Pardo, \emph{Towards a K-theoretic characterization of graded isomorphisms between Leavitt path algebras}. Journal of K-Theory, 14(2), 203-245 (2014). 

\bibitem{A1_6} P. Ara, R. Hazrat, H. Li, A. Sims, \emph{Graded {S}teinberg algebras and their representations}, Algebra Number Theory 12 (1) (2018), 131-172.

\bibitem{A1_4} P. Ara, M.A. Moreno, E. Pardo, Nonstable K-theory for graph algebras, Algebras and Representation Theory 10 (2007), 157-178.

\bibitem{Asmussen} S. R. Asmussen, \emph{{M}arkov Chains}. Applied Probability and Queues. Stochastic Modelling and Applied Probability. Vol. 51, (2003), 3-8.



\bibitem{BS22} W. Bock and A. Sebandal. \emph{A talented monoid view on {L}ie bracket algebras over Leavitt path algebras}. Journal of Algebra and Its Applications,(2022), 2350170.

	
	
\bibitem{A2} L.G. Cordeiro, D. Gon\c{c}alvez, R. Hazrat, \emph{The talented monoid of a directed graph with applications to graph algebras}, preprint
	
	

\bibitem{monoid.q} Y. Give'on, \emph{Normal monoids and factor monoids of commutative monoids}, The University of Michigan, (1963).


\bibitem{LieAlg1} B.C. Hall, \emph{Lie groups, Lie algebras, and representations.} In Quantum Theory for Mathematicians (pp. 333-366). Springer, New York, NY. (2013)


\bibitem{A1} R. Hazrat, H. Li, \emph{The talented monoid of a Leavitt path algebra}, Journal of Algebra 547.  430-455. (2020).


\bibitem{A1_8} R. Hazrat, \emph{The graded Grothendieck group and the classification of Leavitt path algebras}, Mathematische Annalen. 355, 273-325. (2013).


\bibitem{hazratGK} R. Hazrat, A.N. Sebandal, J.P. Vilela, \emph{Graphs with disjoint cycles classification via the talented monoid}, Journal of Algebra. 593.  319-340. (2022).

\bibitem{A3_8} R. Hazrat, L. Vas, \emph{Comparability in the graph monoid}, New York J. of Math, 26 (2020), 1375-1421.


\bibitem{LieAlg2} J. E. Humphreys,  \emph{Introduction to Lie algebras and representation theory} (Vol. 9). Springer Science \& Business Media. (2012).


\bibitem{GKdimension} G.R. Krause, T.H.  Lenagan, Growth of Algebras and Gelfand-Kinillow Dimension, revised edition. Graduate studies in Mathematics, 22, AMS, Providence, RI, 2000.


\bibitem{leavitt} W.G. Leavitt, \emph{
The module type of a ring}.
 Transactions of the American Mathematical Society, 103, pp. 113-130. (1962).
 
 \bibitem{A1_1_7} W.G. Leavitt, The module type of homomorphic images, Duke Math. J. 32, 305-311, (1965).


\bibitem{namzerui} T.G. Nam, Z. Zhang, \emph{Lie solvable Lie Algebras}, Journal of Algebra and Its Applications. (2021).

\bibitem{Nam}  T.G.,Nam, \emph{Simple Lie algebras arising from Steinberg algebras of Hausdorff ample groupoids}. Journal of Algebra, 595, pp.194-215. (2022)

\bibitem{Nelson1977} Nelson, P. R. \emph{Single-shelf library-type Markov chains with infinitely many books.} Journal of Applied Probability, 14(2), 298-308. (1977)

\bibitem{sebandalvilela} A.N. Sebandal, J.P. Vilela, \emph{The Jordan-Hölder Theorem for monoids with group action}, Journal of Algebra and its Applications, (2022).


\bibitem{Shestakov} I. Shestakov, and Z. Zhang,  \emph{Solvability and nilpotency of Novikov algebras.} Communications in Algebra, 48(12), 5412-5420, (2020).

\bibitem{Vas} L. Vas, \emph{Every graded ideal of a Leavitt path algebra is graded isomorphic to a Leavitt path algebra.} Bulletin of the Australian Mathematical Society, 1-9 (2021).

\bibitem{Vas2} L. Vas, \emph{Graded irreducible representations of Leavitt path algebras: a new type and complete classification.} arXiv preprint arXiv:2201.02446, (2022). 

\bibitem{Wagoner} J. B. Wagoner, \emph{Topological Markov chains, $C^*$-algebras, and $K_2$}. Advances in Mathematics, 71(2), 133-185.
Chicago, (1988). 

\end{thebibliography}
\end{document}